\documentclass[11pt,reqno]{amsart}
\usepackage{graphicx,amssymb,mathrsfs,amsmath,amsfonts,appendix}
\usepackage{subfigure}
\usepackage{amsthm}
\usepackage{paralist}
\usepackage{enumitem}
\usepackage[utf8]{inputenc}
\usepackage{comment}
\usepackage{amsmath}
\usepackage{mathtools}
\usepackage{esint}
\usepackage[colorlinks=true, pdfstartview=FitV, linkcolor=blue, citecolor=blue, urlcolor=blue]{hyperref}
\numberwithin{equation}{section}

\newtheorem{theorem}{Theorem}[section]

\newtheorem{remark}{Remark}[section]

\def\erre{\mathbb{R}}

\usepackage[colorinlistoftodos]{todonotes}

\newcommand{\lara}[1]{\todo[inline,color=red!40]{Lara: #1}}
\newcommand{\elisa}[1]{\todo[inline,color=blue!40]{Elisa: #1}}
\newcommand{\luca}[1]{\todo[inline,color=orange!40]{Luca: #1}}
\newcommand{\chris}[1]{\todo[inline,color=green!40]{Chris: #1}}

\newcommand{\xd}{\textrm{d}}

\newcommand{\ep}[0]{\varepsilon}

\newcommand{\embed}{\hookrightarrow}
\def\wstarto{\stackrel{*}{\rightharpoonup}}

\allowdisplaybreaks

\title[Local asymptotics for the nonlocal Swift-Hohenberg equation]
{Local asymptotics for the nonlocal Swift-Hohenberg equation}

\author[Elisa Davoli]{Elisa Davoli}
\address{Institute for Analysis and Scientific Computing, Vienna University of Technology, Wiedner Hauptstrasse 8-10, 1040 Vienna, Austria}
\email{elisa.davoli@tuwien.ac.at}
\author[Christian Kuehn]{Christian Kuehn}
\address{Technical University of Munich (TUM)
Faculty of Mathematics, Boltzmannstr. 3,
85748 Garching bei München, Germany}
\email{ckuehn@ma.tum.de}
\author[Luca Scarpa]{Luca Scarpa}
\address{Department of Mathematics, Politecnico di Milano, Via E. Bonardi 9, 20133 Milano, Italy}
\email{luca.scarpa@polimi.it}
\author[Lara Trussardi]{Lara Trussardi}
\address{Department of Mathematics and Scientific Computing, University of Graz, Heinrichstraße 36, 8010 Graz (Austria)}
\email{lara.trussardi@uni-graz.at}

\begin{document}

\begin{abstract}
The nonlocal-to-local asymptotics investigation for evolutionary problems is a central topic both in the theory of PDEs and in functional analysis. More recently, it became the main core of the mathematical analysis of phase-separation models. In this paper we focus on the Swift-Hohenberg equations which are key benchmark models in pattern formation problems and amplitude equations. We prove well-posedness of the nonlocal Swift-Hohenberg equation, and study the nonlocal-to-local asymptotics with one and two nonlocal contributions under homogeneous Neumann boundary conditions using suitable energy estimates on the nonlocal problems. 
\end{abstract}

\maketitle

\tableofcontents

%%%%%%%%%%%%%%%%%%%%%%%%%%%%%%%%%%%%%%%%%%%%%%%%%%%%%%%%%
\section{Introduction}
\label{sec:intro}

We are interested in the Swift-Hohenberg equation~\cite{CH93, SH77}. In dimension one, the local Swift-Hohenberg equation is a partial differential equation (PDE) for $u(x,t)\in\mathbb{R}$, with $x\in\Omega\subseteq\mathbb{R}$ and $t\geq0$, given by
\begin{equation}
 \label{eq:SHL}
 \partial_t u = [r-(1+\partial_x^2)^2] u + N(u).
\end{equation}
Here, $r$ is a real parameter, the nonlinearity is typically either $N(u)=bu^2-u^3$ or $N(u)=su^3-u^5$ (with parameters $b,s>0$), and $(1+\partial_x^2)^2 u$ 
is a short-hand notation for $u+\partial_x^4 u+2\partial_x^2 u$. The Swift-Hohenberg PDE is widely used in several fields ranging from fluid mechanics, laser physics, chemistry, and biological systems to vast classes of other pattern formation problems~\cite{CH93,LegaMoloneyNewell,KozyreffTlidi}. In the dynamical systems analysis of PDEs, the Swift-Hohenberg PDE has become a standard benchmark model, particularly for the derivation of amplitude equations~\cite{ColletEckmann1,KirrmannSchneiderMielke,Hoyle}. The study of pattern formation instabilities is usually carried out varying the bifurcation parameter $r$. In particular, the solution $u(x,t)\equiv0$ is linearly stable for $r<0$, and for $r=0$ it becomes unstable, leading to the formation of a family of non-homogeneous patterned states. Furthermore, it has been proven that it can exhibit various complicated wave patterns~\cite{Mielke4,Schneider5,LloydSandstedeAvitabileChampneys}. In addition, it has also been discovered that the Swift-Hohenberg equation can support various spatially-localized patterns lying on interlaced parameter space curves referred to homoclinic snaking~\cite{BurkeKnobloch,MD14}. To understand the spatial localization and parameter space structure of the Swift-Hohenberg equation, it is natural to unfold the model including nonlocal terms. The first groundbreaking work in this direction proposing a nonlocal equation has been \cite{MD14}, it is given by 
\begin{equation}
 \label{eq:SHNL}
 \partial_t u = [r-(1+\partial_x^2)]^2 u +N(u)-\gamma  u(x,t)\int_\Omega K(x-y)u(y,t)^2\xd y
\end{equation}
where $K(x)$ is a bounded and symmetric function defined in $\Omega$ with compact support. The parameter $\gamma$ is the coefficient of the nonlocal term. This nonlocal variant has generated significant interest, e.g., leading to recent proofs that locally near the first bifurcation point many nonlocal Swift-Hohenberg equations actually lead to local amplitude equations~\cite{KuehnThrom,KuehnThrom3}. This naturally raised the question regarding a global view, i.e., under which conditions solutions of nonlocal Swift-Hohenberg equations converge to local solutions in the singular limit of the nonlocality converging to a local term. In this paper we study this question and provide concrete answers for a wide variety of nonlocal kernels.

To carry out the analysis, we first note that in both cases (local and nonlocal), the Swift-Hohenberg equation is the gradient flow associated to a suitable energy functional  (see~\cite{MD14}). For \eqref{eq:SHL}, the free energy is given by
\[
 E_L(u) = \int_\Omega \Bigl( \frac{1}{2}u_{xx}^2 -u_x^2 + \frac{1-r}{2}u^2- \int^u N(v)\xd v\Bigr)\xd x,
\]
where $\int^u N(v)\xd v$ denotes a suitable primitive of $N$. We observe that $u_t= -\delta E_L/\delta u$ guarantees convergence towards the equilibrium state as long as $ E_L $ is bounded from below. For~\eqref{eq:SHNL} the energy functional is given by
\[
 E_{NL} = \int_\Omega \Bigl( \frac{1}{2}u_{xx}^2 -u_x^2 + \frac{1-r}{2}u^2- \int^u N(v)\xd v\Bigr)\xd x 
 +\frac{\gamma}{4}\int_\Omega\int_\Omega K(x-y)u^2(x)u^2(y)\xd x\xd y.
\]
Also in this case, one can ensure that the time evolution $ u_t = -\delta E_{NL}/\delta u $ is well posed if $ E_{NL} $ is bounded from below and therefore, in the case with the choice $ N(u) = bu^2-u^3 $, we are going to consider the case $ \gamma> -1 $.

Our goal is to prove the convergence of the solutions of the nonlocal problem~\eqref{eq:SHNL} to the solutions of the local problem~\eqref{eq:SHL} under a wide generality on the setting.
In particular, we work in any space dimension (not just in dimension 1) under typical Neumann boundary conditions
and we handle reasonable classes of the kernel $K$. The idea will be to consider a family of convolution kernels, 
depending on a parameter $\ep$, that peaks around a Dirac delta 
as the $\ep$ vanishes. The key technique
will be to employ suitable energy estimates on the nonlocal problems, which are uniform with respect to $\ep$, in order to guarantee convergence as $\ep\to0$ to the local problem.

The nonlocal-to-local asymptotics investigation for evolutionary 
problems is a central topic both in the theory of PDEs and 
even in functional analysis. It dates back to the 
nonlocal-to-local studies on energy functionals 
carried out by J.~Bourgain, H.~Brezis, 
P.~Mironescu in \cite{BBM, BBM2} 
in relation to Sobolev space theory,
and by Ponce in \cite{ponce04, ponce}
in the context of 
nonlocal Poincar\'e inequalities.
The investigation of nonlocal-to-local asymptotics
in relation to Gamma convergence for 
evolutionary problems 
has been pioneered by 
E.~Sandier and 
S.~Serfaty in \cite{sand-serf}.
More recently, nonlocal-to-local asymptotics 
has been investigated
by some of the authors 
in the context of phase-separation models 
of Cahn-Hilliard type. The pioneering works 
\cite{MRT18, DRST, DST, DST2} initiated such study,
which has now become one of the main cores of 
the mathematical analysis of phase-separation models, 
see \cite{ES, AT} and the references therein.

The paper is organised in two main sections.
Section~\ref{sec:1} presents the study on the Swift-Hohenberg equation with one nonlocal term, while 
Section~\ref{sec:2} deals with the general case of two nonlocal terms.

%%%%%%%%%%%%%%%%%%%%%%%%%%%%%%%%%%%%%%%%%%
\section{Swift-Hohenberg with one nonlocal contribution}
\label{sec:1}
\subsection{Setting of the problem}
\label{sec:intro1}
We consider a family of nonlocal Swift-Hohenberg equations,
depending on a parameter $\ep>0$, on 
a bounded domain $\Omega\subset\mathbb{R}^d$ for $d\geq1$
with homogeneous Neumann boundary conditions of the form 
\begin{align}
 \label{eq:SHNL_d1}
 \partial_t u_\ep + (I+\Delta)^2u_\ep
 = r u_\ep +N(u_\ep)-\gamma  u_\ep K_\ep*u^2_\ep
 \qquad&\text{in } (0,T)\times\Omega\,,\\
 \label{eq:SHNL_d2}
 \partial_{\bf n}u_\ep=\partial_{\bf n}\Delta u_\ep=0 \qquad&\text{on } (0,T)\times\partial\Omega\,,\\
 \label{eq:SHNL_d3}
 u(0)=u_{0,\ep} \qquad&\text{in } \Omega\,.
\end{align}
We are interested in studying the asymptotic behaviour 
of the system \eqref{eq:SHNL_d1}--\eqref{eq:SHNL_d3}
when $\ep\searrow0$, according to different 
possible choices of the family of convolution kernels 
$(K_\ep)_\ep$ and of their scaling with respect to $\ep$.

%%%%%%%%%%%%%%%%%%%%%%%%%%%%%%%%%%%%%%%%%%
\subsection{Assumptions and main results}
\label{sec:main}
We fix here the main assumptions of the work and 
introduce the setting that we will use. Throughout the paper, $\Omega$ is a 
bounded Lipschitz domain in $\mathbb{R}^d$, with $d\in\{1,2,3\}$,
and $T>0$ is a fixed final time. We use the notation 
\begin{align*}
  H^2_{\bf n}(\Omega)&:=
  \left\{ \varphi \in H^2(\Omega):
  \;\partial_{\bf n}\varphi = 0\;
  \text{ a.e.~on } \partial\Omega
  \right\}\,,\\
  H^4_{\bf n}(\Omega)&:=
  \left\{ \varphi \in 
  H^2_{\bf n}(\Omega)\cap H^4(\Omega):
  \;\partial_{\bf n}\Delta\varphi = 0\;
  \text{ a.e.~on } \partial\Omega
  \right\}\,.
\end{align*}
The reason why we introduce such spaces is readily clear.
Indeed, the Laplace operator 
$\Delta$ with homogeneous Neumann 
boundary condition is an unbounded linear operator in 
$L^2(\Omega)$: the spaces
$H^2_{\bf n}(\Omega)$ and $H^4_{\bf n}(\Omega)$ 
coincide exactly 
with the effective domains of $\Delta$ and $\Delta^2$
on $L^2(\Omega)$, respectively. Let us 
also recall that the differential linear operator 
entering the equation \eqref{eq:SHNL_d1} reads
\[
  (I+\Delta)^2:H^4_{\bf n}(\Omega) \to L^2(\Omega)\,,
  \qquad
  (I+\Delta)^2\varphi:=
  \varphi + \Delta^2\varphi + 2\Delta\varphi\,,
  \quad\varphi\in H^4_{\bf n}(\Omega)\,.
\]

We assume the following.
\begin{description}
  \item[A0] $\gamma\in\erre$ and $r\in\erre$
  \item[A1] $N:\mathbb{R}\to\mathbb{R}$ is 
  of class $C^1$
  and we set $\widehat N:\erre\to\erre$ as
  \[
  \widehat N(x):=-\int_0^xN(s)\,\xd s\,, \quad x\in\erre\,.
  \]
  We assume that 
  there exist constants $c_N,C_N>0$ and $p\geq4$ such that
  \begin{align*}
  c_N\frac{|x|^p}p - C_N \leq
  \widehat N(x)
  \qquad\forall\,x\in\mathbb{R}\,,
  \end{align*}
  with the further requirement that
  \[
  \gamma\geq-c_N 
  \qquad\text{if } p=4\,.
  \]
  These conditions are satisfied by 
  several
  choices of the potential $N$: 
  for example, 
  every relevant
  polynomial potential 
  is included.
  Let us stress that 
  this assumption also allows to 
  include possibly superpolynomial
  potentials, 
  as the
  exponential ones for instance.
  \item[A2] $K\in L^1(\erre^d)$ is even, nonnegative, and satisfies
  \[
  \int_{\erre^d}K(y)\,\xd y = 1\,.
  \]
  For every $\ep>0$, we set 
  \[
  K_\ep\in L^1(\erre^d)\,, \qquad
  K_\ep(x):=\frac1{\ep^d}K(x/\ep)\,,
  \quad\text{a.e.~}x\in\erre^d\,.
  \]
  For every $\varphi\in L^1(\Omega)$, we use the notation
  \[
  K_\ep*\varphi\in L^1(\Omega)\,, \qquad
  (K_\ep*\varphi)(x):=
  \int_\Omega K_\ep(x-y)\varphi(y)\,\xd y\,,
  \quad\text{a.e.~}x\in\Omega\,.
  \]
  It is well known that
  \begin{align*}
      \|K_\ep*\varphi\|_{L^p(\Omega)}\leq
      \|\varphi\|_{L^p(\Omega)} \qquad&\forall\,
      \varphi\in L^p(\Omega)\,,
      \quad\forall\,p\in[1,+\infty]\\
      K_\ep*\varphi \to \varphi
      \quad\text{in } L^p(\Omega) \qquad&\forall\,
      \varphi\in L^p(\Omega)\,,
      \quad\forall\,p\in[1,+\infty)\,.
  \end{align*}
  We define the nonlocal 
  component 
  $E_\ep:L^4(\Omega)\to\mathbb{R}$
  of the energy as
  \[
  E_\ep(\varphi):=
  \frac\gamma4
    \int_{\Omega\times\Omega}
    K_\ep(x-y)|\varphi(x)|^2
    |\varphi(y)|^2\,\xd x\,\xd y\,,
  \quad\varphi\in L^4(\Omega)\,,
  \]
  and note that actually $E_\ep$
  is well-defined on $L^4(\Omega)$
  thanks to the H\"older inequality 
  since
  \[
  E_\ep(\varphi)\leq
  \frac\gamma4
  \|K_\ep*u_\ep^2\|_{L^2(\Omega)}
  \|u_\ep^2\|_{L^2(\Omega)}
  \leq
  \frac\gamma4
  \|u_\ep\|_{L^4(\Omega)}^4\,.
  \]
  Let us stress that this inequality,
  together with the assumption on
  $\gamma$ in {\bf A1}, ensure that
  the energy 
  \[
  \varphi\mapsto \int_\Omega
  \widehat N(\varphi(x))\,\xd x
  +E_\ep(\varphi)\,,
  \qquad \varphi\in L^4(\Omega)\,,
  \]
  is always bounded from below
  uniformly in $\ep$
  (also when  $p=4$).
  \item[A3] $u_0\in H^2_{\bf n}(\Omega)$,
  $\widehat N(u_0)\in L^1(\Omega)$,
  and
  $(u_{0,\ep})_\ep\subset H^2_{\bf n}(\Omega)$ are such that,
  for some $M_0>0$,
  \[
  u_{0,\ep} \rightharpoonup u_0 
  \quad\text{in } H^2_{\bf n}(\Omega)
  \quad\text{as } \ep\searrow0\,,
  \qquad
  \|\widehat N(u_{0,\ep})\|_{L^1(\Omega)}\leq 
  M_0 \quad\forall\,\ep>0\,.
  \]
\end{description}

\begin{theorem}
  \label{thm:WP_nl}
  For every $\ep>0$ there exists a unique 
  \[
  u_\ep \in H^1(0,T; L^2(\Omega))
  \cap L^\infty(0,T; H^2_{\bf n}(\Omega)) \cap
  L^2(0,T; H^4_{\bf n}(\Omega))
  \]
  such that $u_\ep(0)=u_{0,\ep}$ and 
  \[
  \partial_t u_\ep + (I+\Delta)^2u_\ep 
  = r u_\ep 
  +N(u_\ep)
  -\gamma  u_\ep K_\ep*u^2_\ep
  \qquad\text{a.e.~in } (0,T)\times\Omega\,.
  \]
\end{theorem}

\begin{theorem}
  \label{thm:conv}
  There exists a unique 
  \[
  u \in H^1(0,T; L^2(\Omega))
  \cap L^\infty(0,T; H^2_{\bf n}(\Omega)) \cap
  L^2(0,T; H^4_{\bf n}(\Omega))
  \]
  such that $u(0)=u_{0}$ and
  \[
  \partial_t u + (I+\Delta)^2u 
  = r u +N(u) -\gamma u^3
  \qquad\text{a.e.~in } (0,T)\times\Omega\,.
  \]
  Moreover, as $\ep\searrow0$ it holds that
  \begin{align*}
      u_\ep \to u \qquad&\text{in }
      C^0([0,T]; H^s(\Omega))
      \quad\forall\,s<2\,,\\
      u_\ep \rightharpoonup u \qquad&\text{in }
      H^1(0,T; L^2(\Omega))
      \cap L^2(0,T; H^4_{\bf n}(\Omega))\,,\\
      u_\ep \wstarto u \qquad&\text{in }
      L^\infty(0,T; H^2_{\bf n}(\Omega))\,.
  \end{align*}
  In particular, it holds that 
  \begin{align*}
  u_\ep \to u \quad&\text{in }
      C^0([0,T]; C^{1,\alpha}(\overline\Omega))
      \quad\forall\,\alpha
      \in(0,1/2)
      &&\text{if $d=1$}\,,\\
  u_\ep \to u \quad&\text{in }
      C^0([0,T]; C^{0,\alpha}(\overline\Omega))
      \quad\forall\,\alpha
      \in(0,1)
      &&\text{if $d=2$}\,,\\
  u_\ep \to u \quad&\text{in }
      C^0([0,T]; C^{0,\alpha}(\overline\Omega))
      \quad\forall\,\alpha
      \in(0,1/2)
      &&\text{if $d=3$}\,.
  \end{align*}
\end{theorem}

Before continuing with the proofs of our main results, let us briefly comment on our assumptions {\bf A0}-{\bf A3}. The assumptions {\bf A0}, {\bf A1} and {\bf A3} are natural from the viewpoint of the classical Swift-Hohenberg regarding the structure of the nonlinearity and regularity conditions. The nonlocal kernel $K_\varepsilon$ in {\bf A2} has to be an approximate identity as $\varepsilon\rightarrow 0$ to study the nonlocal-to-local transition. Furthermore, a kernel in the quadratic nonlinearity has recently not only appeared in the context of the nonlocal Swift-Hohenberg equation~\cite{KuehnThrom3} but also features in other nonlocal reaction-diffusion equations~\cite{Becketal2} including most prominently the nonlocal Fisher-KPP equation~\cite{BerestyckiNadinPerthameRyzhik,AchleitnerKuehn,NadinPerthameTang,KuehnTkachov}. 

%%%%%%%%%%%%%%%%%%%%%%%%%%%%%%%%%%%%%%%%%%
\subsection{Well-posedness of the nonlocal problem and nonlocal-to-local convergence}
\label{sec:WP_NL}

The proof of existence of a solution to the nonlocal problem are based on the uniform estimates done in~\ref{s:un_est_1}, with $\ep>0$ fixed.
In order to obtain the well-posedness, one can follow the same techniques used in~\cite{Gio16}, using an approximation scheme (of Galerkin type).

%\luca{Si potrebbe dire che le stime che facciamo nella sezione 4.1 permettono di otterenere well-posedness del nonlocale se replicate  su uno schema di approssimazione (tipo Galerkin). Cosi si potreebbero fondere sezione 3 e 4.}

%\subsection{Approximation}

%\subsection{Uniform estimates}

%\subsection{Passage to the limit as $\lambda\searrow0$}

%%%%%%%%%%%%%%%%%%%%%%%%%%%%%%%%%%%%%%%%%%
%\section{Nonlocal-to-local convergence with one nonlocal contribution}
%\label{sec:conv}

\subsubsection{Uniform estimates}\label{s:un_est_1}
Testing \eqref{eq:SHNL_d1} by $\partial_t u_\ep$ we get
\begin{align*}
    &\int_0^t\|\partial_t
    u_\ep(s)\|^2_{L^2(\Omega)}\, \xd s
    +\frac12\|u_\ep(t)\|_{L^2(\Omega)}^2
    +\frac12\|\Delta u_\ep(t)\|_{L^2(\Omega)}^2
    +\int_\Omega
    \widehat N(u_\ep(t,x))
    \,\xd x
    +E_\ep(u_\ep(t))\\
    &=\frac12\|u_{0,\ep}\|_{L^2(\Omega)}^2
    +\frac12\|\Delta u_{0,\ep}\|_{L^2(\Omega)}^2
    +\int_\Omega
    \widehat N(u_{0,\ep}(x))
    \,\xd x
    +E_\ep(u_{0,\ep})\\
    &\qquad
    -2\int_0^t\int_\Omega
    \Delta u_\ep(s,x)
    \partial_t u_\ep(s,x)\,\xd x\,\xd s
    +r\int_0^t\int_\Omega
    u_\ep(s,x)\partial_t u_\ep(s,x)\,\xd x\,\xd s\,.
\end{align*}
Now, the first three terms on the
right-hand side are
uniformly bounded in $\ep$ thanks to 
assumption {\bf A3}.
Moreover, recalling {\bf A2} 
we have 
\[
  E_\ep(u_{0,\ep})\leq
  \|u_{0,\ep}\|_{L^4(\Omega)}^4\,,
\]
which is uniformly bounded in 
$\ep$ again thanks to {\bf A3}.
Furthermore, by the weighted Young
inequality it holds
that 
\begin{align*}
    &-2\int_0^t\int_\Omega
    \Delta u_\ep(s,x)
    \partial_t u_\ep(s,x)\,\xd x\,\xd s
    +r\int_0^t\int_\Omega
    u_\ep(s,x)\partial_t u_\ep(s,x)\,\xd x\,\xd s\\
    &\leq \frac12
    \int_0^t\|\partial_t
    u_\ep(s)\|^2_{L^2(\Omega)}\, \xd s
    +4\int_0^t\|\Delta
    u_\ep(s)\|^2_{L^2(\Omega)}\, \xd s
    +r^2\int_0^t\|
    u_\ep(s)\|^2_{L^2(\Omega)}\, \xd s\,.
\end{align*}
Putting all this information 
together and
possibly updating the value of the constant $M$,
independently of $\ep$, 
we are left with 
\begin{align*}
    &\frac12\int_0^t\|\partial_t
    u_\ep(s)\|^2_{L^2(\Omega)}\, \xd s
    +\frac12\|u_\ep(t)\|_{L^2(\Omega)}^2
    +\frac12\|\Delta u_\ep(t)\|_{L^2(\Omega)}^2
    +\int_\Omega\widehat N(u_\ep(t,x))\,\xd x
    +E_\ep(u_\ep(t))\\
    &\leq
    M
    \left(1+
    \int_0^t\|\Delta
    u_\ep(s)\|^2_{L^2(\Omega)}\, \xd s
    +\int_0^t\|
    u_\ep(s)\|^2_{L^2(\Omega)}\, \xd s
    \right)\,.
\end{align*}
Recalling that as a consequence 
of {\bf A1--A2} the energy
contribution 
\[
  \int_\Omega\widehat N(u_\ep(t,x))\,\xd x
    +E_\ep(u_\ep(t))
\]
is bounded from below uniformly in $\ep$, 
the Gronwall lemma and elliptic regularity yield then
\begin{align}
    \label{est1}
    \|u_\ep\|_{H^1(0,T; L^2(\Omega))\cap 
    C^0([0,T]; H^2_{\bf n}(\Omega))} &\leq M\,.
\end{align}
At this point, 
the H\"older inequality and the 
continuous inclusion 
$H^2_{\bf n}(\Omega)\hookrightarrow
C^0(\overline\Omega)$
yield
\[
  \|u_\ep 
  K_\ep*u_\ep^2\|_{C^0(\overline\Omega)}
  \leq
  \|u_\ep\|_{C^0(\overline\Omega)}
  \|K_\ep*u_\ep^2\|_{C^0(\overline\Omega)}
  \leq 
  \|u_\ep\|_{C^0(\overline\Omega)}^3
  \leq M
  \|u_\ep\|_{H^2_{\bf n}(\Omega)}^3
  \,,
\]
which in turn implies by \eqref{est1} that
\begin{equation}
    \label{est2}
    \|u_\ep 
  K_\ep*u_\ep^2 \|_{C^0([0,T];
  C^0(\overline\Omega))} \leq M\,.
\end{equation}
Furthermore, 
since $H^2_{\bf n}(\Omega)\embed
C^0(\overline\Omega)$
and $N\in C^1(\erre)$, the estimate
\eqref{est1} yields also
\begin{equation}
    \label{est3}
    \|N(u_\ep)\|_{C^0([0,T]; C^0(\overline\Omega))}
    +\|N'(u_\ep)\|_{C^0([0,T]; C^0(\overline\Omega))}
    \leq M\,.
\end{equation}
Eventually,
by comparison in the equation \eqref{eq:SHNL_d1},
we infer that
\[
  \|\Delta^2u_\ep\|_{
    L^2(0,T;L^2(\Omega))} \leq M\,,
\]
so that elliptic regularity and the boundary conditions in \eqref{eq:SHNL_d2} yield
also
\begin{equation}
\label{est4}
    \|u_\ep\|_{
    L^2(0,T;H^4_{\bf n}(\Omega))} \leq M\,.
\end{equation}

\subsubsection{Passage to the limit}
From the estimates \eqref{est1}--\eqref{est4}
and the classical Aubin-Lions-Simon
compactness results, we 
infer that there exists
\[
  u\in 
  H^1(0,T; L^2(\Omega))
  \cap L^\infty(0,T; H^2_{\bf n}(\Omega)) \cap
  L^2(0,T; H^4_{\bf n}(\Omega))
\]
such that as $\ep\searrow0$
(possibly on a non-relabelled subsequence)
\begin{align*}
      u_\ep \to u \qquad&\text{in }
      C^0([0,T]; H^s(\Omega))
      \quad\forall\,s<2\,,\\
      u_\ep \rightharpoonup u \qquad&\text{in }
      H^1(0,T; L^2(\Omega))
      \cap L^2(0,T; H^4_{\bf n}(\Omega))\,,\\
      u_\ep \wstarto u \qquad&\text{in }
      L^\infty(0,T; H^2_{\bf n}(\Omega))\,.
\end{align*}
As $H^2_{\bf n}(\Omega)\embed C^0(\overline \Omega)$, 
by the continuity of $N$
this implies that 
\[
  N(u_\ep(t,x))\to N(u(t,x))
  \qquad\forall\,(t,x)
  \in[0,T]\times\Omega\,,
\]
hence from \eqref{est3} we deduce 
\[
  N(u_\ep)\wstarto N(u)
  \qquad\text{in } L^\infty((0,T)\times \Omega)\,.
\]
We are only left to show how to pass to
the limit in the nonlocal term.
To this end, first of all
we have, for all $q\in[1,+\infty)$
\begin{align*}
  \|K_\ep*u_\ep^2 - u^2\|_{L^q((0,T)\times\Omega)}&\leq
  \|K_\ep*(u_\ep^2-u^2)\|_{L^q((0,T)\times\Omega)}
  +\|K_\ep*u^2-u^2\|_{L^q((0,T)\times\Omega)}\\
  &\leq
  \|u_\ep^2-u^2\|_{L^q((0,T)\times\Omega)}
  +\|K_\ep*u^2-u^2\|_{L^q((0,T)\times\Omega)}\to 0\,,
\end{align*}
so that 
\[
  K_\ep*u_\ep^2 \to u^2
  \qquad\text{in } L^q((0,T)\times\Omega)
  \quad\forall\,q\in[1,+\infty)\,.
\]
Recalling the convergence of $(u_\ep)_\ep$, this shows in particular that
\[
  -\gamma u_\ep K_\ep*u_\ep^2 \to 
  -\gamma u^3 \qquad\text{in }
  L^2(0,T; L^2(\Omega))\,.
\]
Letting then $\ep\searrow0$
we deduce that $u$ solves the local 
equation in Theorem~\ref{thm:conv}.
As the local equation has a unique 
solution, we infer that
the convergences hold along the entire
sequence, and the proof of Theorem~\ref{thm:conv}
is concluded.

The uniqueness at the limit is guaranteed by the result of~\cite{Gio16}, assuming $N$ of class $C^3$.
Moreover, since we are working with Neumann boundary conditions, instead of Dirichlet, we should additionally have $u$ with zero-mean in order to guarantee uniqueness.

\section{Swift-Hohenberg with two nonlocal contributions}
\label{sec:2}
%\section{Nonlocal-to-local convergence with two nonlocal contributions}
%%%%%%%%%%%%%%%%%%%%%%%%%%%%%%%%%%%%%%%%%%%%%%%%%%%%%%%%%%%%%%%%%%%%%%%%%%%%%%%
\subsection{Setting of the problem}
We consider a family of nonlocal Swift-Hohenberg equations,
depending on a parameter $\ep>0$, on 
a bounded domain $\Omega\subset\mathbb{R}^d$ for $d\geq1$
with Neumann boundary conditions in the form 
\begin{align}
 \label{eq:SHNL2_d1}
 \partial_t u_\ep + (I+\Delta)^2u_\ep
 = r u_\ep -u_\ep (Q_\ep*u_\ep^p)+  u_\ep (K_\ep*u^q_\ep)
 \qquad&\text{in } (0,T)\times\Omega\,,\\
 \label{eq:SHNL2_d2}
 \partial_{\bf n}u_\ep=\partial_{\bf n}\Delta u_\ep=0 \qquad&\text{on } (0,T)\times\partial\Omega\,,\\
 \label{eq:SHNL2_d3}
 u(0)=u_{0,\ep} \qquad&\text{in } \Omega\,,
\end{align}
for $p<q$ and $q$ even. Interesting cases of applications include both the case $p=1$ and $q=2$, leading in the limit to a quadratic-cubic nonlinearity, and the setting $p=2$ and $q=4$, yielding a cubic-quintic term.
We are interested in studying the asymptotic behaviour 
of the system \eqref{eq:SHNL2_d1}--\eqref{eq:SHNL2_d3}
when $\ep\searrow0$, according to different 
possible choices of the families of convolution kernels 
$(K_\ep)_\ep$ and $(Q_\ep)_\ep$, and of their scalings with respect to $\ep$.
\subsection{Assumptions and main results}
We fix here the main assumptions of the work and 
introduce the setting that we will use. We assume the following.
\begin{description}
  \item[H0] $K\in L^1(\erre^d)$ is even, nonnegative, and satisfies
  \[
  \int_{\erre^d}K(y)\,\xd y = 1\,.
  \]
  For every $\ep>0$, we set 
  \[
  K_\ep\in L^1(\erre^d)\,, \qquad
  K_\ep(x):=\frac1{\ep^d}K(x/\ep)\,,
  \quad\text{a.e.~}x\in\erre^d\,.
  \]
  For every $\varphi\in L^1(\Omega)$, we use the notation
  \[
  K_\ep*\varphi\in L^1(\Omega)\,, \qquad
  (K_\ep*\varphi)(x):=
  \int_\Omega K_\ep(x-y)\varphi(y)\,\xd y\,,
  \quad\text{a.e.~}x\in\Omega\,.
  \]
  It is well known that
  \begin{align*}
      \|K_\ep*\varphi\|_{L^s(\Omega)}\leq
      \|\varphi\|_{L^s(\Omega)} \qquad&\forall\,
      \varphi\in L^s(\Omega)\,,
      \quad\forall\,s\in[1,+\infty]\\
      K_\ep*\varphi \to \varphi
      \quad\text{in } L^s(\Omega) \qquad&\forall\,
      \varphi\in L^s(\Omega)\,,
      \quad\forall\,s\in[1,+\infty)\,.
  \end{align*}
  We define the nonlocal 
  component
  $E_\ep^K:L^{q+2}(\Omega)\to\mathbb{R}$
  of the energy as
  \begin{align*}
  E_\ep^K(\varphi)&:=
  \frac14
    \int_{\Omega\times\Omega}
    K_\ep(x-y)|\varphi(x)|^{2}
    |\varphi(y)|^{q}\,\xd x\,\xd y\,,
  \quad\varphi\in L^{q+2}(\Omega)\,.
  \end{align*}
  Arguing as in Section \ref{sec:main}, we see by H\"older's inequality that this energy term is well-defined on $L^{q+2}$.
  \item[H1] $Q\in L^1(\erre^d)$ is even.
%  , nonnegative, satisfies
%  \[
%  \int_{\erre^d}Q(y)\,\xd y = 1\,,
%  \]
%  as well as
%  $$K-Q\geq 0\quad\text{almost everywhere in } \mathbb{R}.$$
%  \chris{Weaken this assumption to avoid pointwise statement but rather make it an integral statement of one kernel winning on average against the other. Maybe just assume that one energy dominates the other...}
 % \elisa{After the discussion on 08.03.2024: eliminate the assumption of $Q$ and $K$ have unitary $L^1$ norm. Replace all assumptions above with: $K$ nonnegative and $|Q|\leq K$. Lemma 5.1 and Section 5.5 need to be amended accordingly}
 % \lara{We assume that one energy dominates the other one. Add a remark saying that if $|Q|\leq K$ then one can prove that one energy dominates the other one. }
 For every $\ep>0$, we set 
  \[
  Q_\ep\in L^1(\erre^d)\,, \qquad
  Q_\ep(x):=\frac1{\ep^d}Q(x/\ep)\,,
  \quad\text{a.e.~}x\in\erre^d\,.
  \]
  For every $\varphi\in L^1(\Omega)$, the notation
  \[
  Q_\ep*\varphi\in L^1(\Omega)
  \]
  is defined analogously to $K_\ep*\varphi$.
  We define the nonlocal 
  component
  $E_\ep^Q:L^{p+2}(\Omega)\to\mathbb{R}$
  of the energy as
  \begin{align*}
  E_\ep^Q(\varphi)&:=
  \frac14
    \int_{\Omega\times\Omega}
    Q_\ep(x-y)|\varphi(x)|^{2}
    |\varphi(y)|^{p}\,\xd x\,\xd y\,,
  \quad\varphi\in L^{p+2}(\Omega)\,.
  \end{align*}
  Once more the well-definitness of this energy term is a direct consequence of H\"older's inequality. 
  We assume that there exist constants $C>0$ and $c\in(0,\frac12)$,
  independent of $\ep$,
  such that
  \begin{equation}
  \label{eq:estimate-eQ}
  |E^Q_\ep(\varphi)|\leq C(1 + E^K_\ep(\varphi)) + c\|\varphi\|_{L^2(\Omega)}^2
  \quad\forall\,\varphi\in L^{q+2}(\Omega)\,,
  \quad\forall\,\ep>0\,.
  \end{equation}
  \item[H2] $u_0\in H^2(\Omega)$,
  $\widehat N(u_0)\in L^1(\Omega)$,
  and
  $(u_{0,\ep})_\ep\subset H^2(\Omega)$ are such that,
  for some $M_0>0$,
  \[
  u_{0,\ep} \rightharpoonup u_0 
  \quad\text{in } H^2(\Omega)
  \quad\text{as } \ep\searrow0\,,
  \qquad
  E_\ep^K(u_{0,\ep})+E_\ep^Q(u_{0,\ep})\leq 
  C_0 \quad\forall\,\ep>0\,.
  \]
\end{description}

As above for {\bf A0}-{\bf A2}, the assumptions {\bf H0}-{\bf H2} are quite natural. The convolution kernel in the cubic nonlinearity has recently appeared in the nonlocal Swift-Hohenberg equation~\cite{MD14,KuehnThrom3}. It can be viewed as analogous to the Fisher-KPP case, just for the nonlocal Ginzburg-Landau (or Allen-Cahn, or Nagumo) equation. Furthermore, it has been proposed in a cubic nonlocal variant of the nonlinear Schr\"odinger equations~\cite{TsilifisKevrekidisRothos,Becketal2}. 

\begin{remark}
Let us point that a sufficient condition for the inequality 
\eqref{eq:estimate-eQ} is that $|Q|\leq K$ almost everywhere.
Indeed, by definition of $Q_\ep$ and $K_\ep$ we have 
\begin{align*}
  |E^Q_\ep(\varphi)|
  &\leq \frac14\int_{\Omega}\int_{\Omega}
  |Q_\ep(x,y)||\varphi(x)|^2|\varphi(y)|^p\,\xd{x}\xd{y}
  \leq \frac14\int_{\Omega}\int_{\Omega}
  K_\ep(x,y)|\varphi(x)|^2|\varphi(y)|^p\,\xd{x}\xd{y}\,.
\end{align*}
Now, since $p<q$, by the weighted Young inequality, for every $\delta>0$
there exists a constant $C_\delta>0$ such that 
\[
|a|^p\leq C_\delta|a|^q + \delta \quad\forall\,a\in\mathbb R\,.
\]
Hence, putting everything together, the required inequality \eqref{eq:estimate-eQ} follows by choosing $\delta<2$.
\end{remark}

We conclude this subsection with the statement of our main results.

\begin{theorem}
  \label{thm:WP2_nl}
  For every $\ep>0$ there exists a unique 
  \[
  u_\ep \in H^1(0,T; L^2(\Omega))
  \cap L^\infty(0,T; H^2_{\bf n}(\Omega)) \cap
  L^2(0,T; H^4_{\bf n}(\Omega))
  \]
  such that $u_\ep(0)=u_{0,\ep}$ and 
  \[
  \partial_t u_\ep + (I+\Delta)^2u_\ep 
  = r u_\ep 
  -u_\ep (Q_\ep*u_\ep^p)+  u_\ep (K_\ep*u^q_\ep)
  \qquad\text{a.e.~in } (0,T)\times\Omega\,.
  \]
\end{theorem}

\begin{theorem}
  \label{thm:conv2}
  There exists a unique 
  \[
  u \in H^1(0,T; L^2(\Omega))
  \cap L^\infty(0,T; H^2_{\bf n}(\Omega)) \cap
  L^2(0,T; H^4_{\bf n}(\Omega))
  \]
  such that $u(0)=u_{0}$ and
  \[
  \partial_t u + (I+\Delta)^2u 
  = r u +u^{p+1}-u^{q+1}
  \qquad\text{a.e.~in } (0,T)\times\Omega\,.
  \]
  Moreover, as $\ep\searrow0$ it holds that
  \begin{align*}
      u_\ep \to u \qquad&\text{in }
      C^0([0,T]; H^s(\Omega))
      \quad\forall\,s<2\,,\\
      u_\ep \rightharpoonup u \qquad&\text{in }
      H^1(0,T; L^2(\Omega))
      \cap L^2(0,T; H^4_{\bf n}(\Omega))\,,\\
      u_\ep \wstarto u \qquad&\text{in }
      L^\infty(0,T; H^2_{\bf n}(\Omega))\,.
  \end{align*}
  In particular, it holds that 
  \begin{align*}
  u_\ep \to u \quad&\text{in }
      C^0([0,T]; C^{1,\alpha}(\overline\Omega))
      \quad\forall\,\alpha
      \in(0,1/2)
      &&\text{if $d=1$}\,,\\
  u_\ep \to u \quad&\text{in }
      C^0([0,T]; C^{0,\alpha}(\overline\Omega))
      \quad\forall\,\alpha
      \in(0,1)
      &&\text{if $d=2$}\,,\\
  u_\ep \to u \quad&\text{in }
      C^0([0,T]; C^{0,\alpha}(\overline\Omega))
      \quad\forall\,\alpha
      \in(0,1/2)
      &&\text{if $d=3$}\,.
  \end{align*}
\end{theorem}

%%%%%%%%%%%%%%%%%%%%%%%%%%%%%%%%%%%%%%%%%%
\subsection{Well-posedness of the nonlocal problem and nonlocal-to-local convergence}

As for the case with one nonlocal contribution, the existence of a solution to the nonlocal problem can be obtained following the same approach as in~\cite{Gio16}. The uniform estimates computed here allow to get the well-posedness.

%\luca{Possiamo fare come all'inizio unendo le due sezioni.}
%\subsection{Approximation}

%\subsection{Uniform estimates}

%\subsection{Passage to the limit as $\lambda\searrow0$}

%%%%%%%%%%%%%%%%%%%%%%%%%%%%%%%%%%%%%%%%%%
%%%%%%%%%%%%%%%%%%%%%%%%%%%%%%%%%%%%%%%%%%%%%%%%%%%%
%\subsection{Nonlocal-to-local convergence with two convolution terms}

\subsubsection{Uniform estimates}
Testing \eqref{eq:SHNL2_d1} by $\partial_t u_\ep$ we get
\begin{align*}
    &\int_0^t\|\partial_t
    u_\ep(s)\|^2_{L^2(\Omega)}\, \xd s
    +\frac12\|u_\ep(t)\|_{L^2(\Omega)}^2
    +\frac12\|\Delta u_\ep(t)\|_{L^2(\Omega)}^2
    +E_\ep^K(u_\ep(t))\\
    &=\frac12\|u_{0,\ep}\|_{L^2(\Omega)}^2
    +\frac12\|\Delta u_{0,\ep}\|_{L^2(\Omega)}^2
    +E_\ep^Q(u_{0,\ep})
    +E_\ep^K(u_{0,\ep})+\left|E_\ep^Q(u_\ep(t))\right|\\
    &\qquad
    -2\int_0^t\int_\Omega
    \Delta u_\ep(s,x)
    \partial_t u_\ep(s,x)\,\xd x\,\xd s
    +r\int_0^t\int_\Omega
    u_\ep(s,x)\partial_t u_\ep(s,x)\,\xd x\,\xd s\,.
\end{align*}
The first four terms on the
right-hand side are
uniformly bounded in $\ep$ thanks to 
assumption {\bf H2}. By the weighted Young
inequality it holds
that 
\begin{align*}
    &-2\int_0^t\int_\Omega
    \Delta u_\ep(s,x)
    \partial_t u_\ep(s,x)\,\xd x\,\xd s
    +r\int_0^t\int_\Omega
    u_\ep(s,x)\partial_t u_\ep(s,x)\,\xd x\,\xd s\\
    &\leq \frac12
    \int_0^t\|\partial_t
    u_\ep(s)\|^2_{L^2(\Omega)}\, \xd s
    +4\int_0^t\|\Delta
    u_\ep(s)\|^2_{L^2(\Omega)}\, \xd s
    +r^2\int_0^t\|
    u_\ep(s)\|^2_{L^2(\Omega)}\, \xd s\,.
\end{align*}

Hypothesis \eqref{eq:estimate-eQ} yields
\begin{equation}
\label{eq:estimate-eQ-applied}
|E_\ep^Q(u_\ep(t))|\leq \frac{|\Omega|}{4}+E_\ep^K(u_\ep(t))+\frac{1}{4}\|u_\ep(t)\|_{L^2(\Omega)}^2.
\end{equation}
Putting all this information 
together and
possibly updating the value of the constant $M$,
independently of $\ep$, 
we are left with 
\begin{align*}
    &\frac12\int_0^t\|\partial_t
    u_\ep(s)\|^2_{L^2(\Omega)}\, \xd s
    +\frac12\|u_\ep(t)\|_{L^2(\Omega)}^2
    +\frac12\|\Delta u_\ep(t)\|_{L^2(\Omega)}^2
    \\
    &\leq
    M
    \left(1+
    \int_0^t\|\Delta
    u_\ep(s)\|^2_{L^2(\Omega)}\, \xd s
    +\int_0^t\|
    u_\ep(s)\|^2_{L^2(\Omega)}\, \xd s
    \right)\,.
\end{align*}
By the Gronwall lemma and by elliptic regularity we infer that
\begin{align}
    \label{est1_2}
    \|u_\ep\|_{H^1(0,T; L^2(\Omega))\cap 
    C^0([0,T]; H^2_{\bf n}(\Omega))} &\leq M\,.
\end{align}
At this point, 
the H\"older inequality and the 
continuous inclusion 
$H^2_{\bf n}(\Omega)\hookrightarrow
C^0(\overline\Omega)$
yield
\[
  \|u_\ep 
  K_\ep*u_\ep^q\|_{C^0(\overline\Omega)}
  \leq
  \|u_\ep\|_{C^0(\overline\Omega)}
  \|K_\ep*u_\ep^q\|_{C^0(\overline\Omega)}
  \leq 
  \|u_\ep\|_{C^0(\overline\Omega)}^{q+1}
  \leq M
  \|u_\ep\|_{H^2_{\bf n}(\Omega)}^{q+1}
  \,,
\]
which in turn implies by \eqref{est1_2} that
\begin{equation}
    \label{est2-two}
    \|u_\ep 
  K_\ep*u_\ep^q \|_{C^0([0,T];
  C^0(\overline\Omega))} \leq M\,.
\end{equation}
Analogously, we obtain that
\begin{equation}
    \label{est2-two-bis}
    \|u_\ep 
  Q_\ep*u_\ep^p \|_{C^0([0,T];
  C^0(\overline\Omega))} \leq M\,.
\end{equation}
Eventually,
by comparison in the equation \eqref{eq:SHNL_d1},
we infer that
\[
  \|\Delta^2u_\ep\|_{
    L^2(0,T;L^2(\Omega))} \leq M\,,
\]
so that elliptic regularity and the boundary conditions in \eqref{eq:SHNL2_d2} yield
also
\begin{equation}
\label{est4-two}
    \|u_\ep\|_{
    L^2(0,T;H^4_{\bf n}(\Omega))} \leq M\,.
\end{equation}

\subsubsection{Passage to the limit}
From the estimates \eqref{est1_2}--\eqref{est4-two}
and the classical Aubin-Lions-Simon
compactness results, we 
infer that there exists
\[
  u\in 
  H^1(0,T; L^2(\Omega))
  \cap L^\infty(0,T; H^2_{\bf n}(\Omega)) \cap
  L^2(0,T; H^4_{\bf n}(\Omega))
\]
such that as $\ep\searrow0$
(possibly on a non-relabelled subsequence)
\begin{align*}
      u_\ep \to u \qquad&\text{in }
      C^0([0,T]; H^s(\Omega))
      \quad\forall\,s<2\,,\\
      u_\ep \rightharpoonup u \qquad&\text{in }
      H^1(0,T; L^2(\Omega))
      \cap L^2(0,T; H^4_{\bf n}(\Omega))\,,\\
      u_\ep \wstarto u \qquad&\text{in }
      L^\infty(0,T; H^2_{\bf n}(\Omega))\,.
\end{align*}
To conclude, it remains only to show how to pass to
the limit in the nonlocal terms.
To this end, first of all
we have, for all $s\in[1,+\infty)$
\begin{align*}
  \|K_\ep*u_\ep^q - u^q\|_{L^s((0,T)\times\Omega)}&\leq
  \|K_\ep*(u_\ep^q-u^q)\|_{L^s((0,T)\times\Omega)}
  +\|K_\ep*u^q-u^q\|_{L^s((0,T)\times\Omega)}\\
  &\leq
  \|u_\ep^q-u^q\|_{L^s((0,T)\times\Omega)}
  +\|K_\ep*u^q-u^q\|_{L^s((0,T)\times\Omega)}\to 0\,,
\end{align*}
so that 
\[
  K_\ep*u_\ep^q \to u^q
  \qquad\text{in } L^s((0,T)\times\Omega)
  \quad\forall\,s\in[1,+\infty)\,.
\]
Recalling the convergence of $(u_\ep)_\ep$, this shows in particular that
\[
  u_\ep K_\ep*u_\ep^q \to 
  u^{q+1} \qquad\text{in }
  L^2(0,T; L^2(\Omega))\,.
\]
Analogously, we obtain that
\[
  u_\ep Q_\ep*u_\ep^p \to 
  u^{p+1} \qquad\text{in }
  L^2(0,T; L^2(\Omega))\,.
\]
Letting then $\ep\searrow0$
we deduce that $u$ solves the local 
equation in Theorem~\ref{thm:conv2}.
As the local equation has a unique 
solution, we infer that
the convergences hold along the entire
sequence, and the proof of Theorem~\ref{thm:conv2}
is concluded.
As previously, the uniqueness at the limit is guaranteed by the result in~\cite{Gio16}, assuming $N$ to be sufficiently regular and by choosing $u$ with zero-mean.

\section*{Acknowledgements}

The research of E.D. has been funded by the Austrian Science Fund (FWF) projects 10.55776/F65, 10.55776/Y1292, 10.55776/P35359, and 10.55776/F100800.
L.S.~is member of Gruppo Nazionale 
per l'Analisi Matematica, la Probabilit\`a 
e le loro Applicazioni (GNAMPA), 
Istituto Nazionale di Alta Matematica (INdAM).
The present research is also part of the activities of 
``Dipartimento di Eccellenza 2023-2027'' of Politecnico di Milano.
For open-access purposes, the authors have applied a CC BY public copyright
license to any author-accepted manuscript version arising from this
submission.

\bibliographystyle{abbrv}
\bibliography{ref}
\end{document}